\documentclass{amsart}

\usepackage[utf8]{inputenc}
\usepackage[T1]{fontenc}
\usepackage{eulervm}
\usepackage{tgpagella}

\usepackage[numbers]{natbib}

\usepackage{amsmath}
\usepackage{amssymb}
\usepackage{amsthm}
\usepackage{amsfonts}
\usepackage{mathrsfs}
\usepackage{xspace}
\usepackage{tikz}
\usepackage{nicefrac}
\usepackage{fixmath}
\usepackage{paralist}
\usepackage{ellipsis}
\usepackage{mathtools}
\usepackage{graphicx,subfigure,epic,eepic}
\usepackage{array}

\usepackage{fixltx2e}
\usepackage[expansion=false,final]{microtype}

\usepackage{MnSymbol}

\newtheorem{thm}{Theorem}
\newtheorem{lem}{Lemma}

\newtheorem{prop}{Proposition}

\DeclareMathOperator{\Pic}{Pic}

\newcommand{\Z}{\mathbb{Z}\xspace}

\newcommand{\C}{\mathbb{C}\xspace}

\newcommand{\bP}{\mathbb{P}\xspace}
\newcommand{\cO}{\mathcal{O}\xspace}
\newcommand{\sD}{\mathscr{D}\xspace}

\title{Strong Franchetta Conjecture for Linear Systems}
\author{Matthew Woolf}
\begin{document}

\begin{abstract}
In this paper, we study rational sections of the relative Picard variety of a linear system on a smooth projective variety. Specifically, we prove that if the linear system is basepoint-free and the locus of non-integral divisors has codimension at least two, then all rational sections of the relative Picard variety come from restrictions of line bundles on the variety.
\end{abstract}

\maketitle

\section{Introduction}

The original Franchetta conjecture was that the only \emph{natural} lines bundles on curves are multiples of the canonical bundle. There are two ways to make this into a precise statement. The weak Franchetta conjecture says that the relative Picard group of the universal curve over $M_g$, i.e.~the group of line bundles on the total space of the universal curve modulo those pulled back from $M_g$, is generated by the relative canonical bundle. The strong Franchetta conjecture says that the only rational sections of the universal Picard varieties $J^d$ come from multiples of the canonical bundle.

Any line bundle on the universal curve over $M_g$ gives rise to a rational section of the universal Picard variety, so the strong Franchetta conjecture implies the weak Franchetta conjecture. The weak Franchetta conjecture follows from the calculation of the Picard group of the universal curve due to Harer and Arbarello-Cornalba. From the weak Franchetta conjecture, it is possible to deduce that given a rational section of the relative Picard variety, some multiple of it is a multiple of the canonical bundle (see proposition \ref{brauertorsion}). Using the weak Franchetta conjecture, there have been a number of proofs of the strong Franchetta conjecture.

The proofs by Mestrano \citep{mestrano} and Kouvidakis \citep{kouvidakis} work by showing that the degree of a rational section must be a multiple of $2g-2$ to reduce to the case of $J^0$, where any section must \emph{a fortiori} be torsion, and then show that the monodromy action on torsion line bundles has no nonzero fixed point.

Given a linear system of curves on a surface, a natural conjecture to make is that the only ``natural line bundles'' are restrictions of line bundles on the surface (at least when the linear system is basepoint-free). We can also make the same conjecture about linear systems in varieties of dimension greater than two. The purpose of this paper is to show that with one additional necessary hypothesis, this conjecture is correct.

\begin{thm}\label{mainresult}
Let $X$ be a smooth projective variety over $\C$, $|D|$ a basepoint-free linear system of divisors on $X$ such that the locus of divisors which are either reducible or non-reduced has codimension at least two. Then rational sections of the relative Picard variety all come from restricting line bundles on $X$.
\end{thm}

This theorem can be thought of as a weakening of the Lefschetz hyperplane theorem -- there are no hypotheses about ampleness or dimension, but we only get a result about the relative Picard group, not the Picard groups of each divisor. In the case of curves, though, where the Lefschetz hyperplane theorem tells you very little, this result provides the most information. For example, by considering curves in K3 surfaces which generate the Picard group, we see immediately that the degree of any rational section of the universal Picard variety over the moduli space of curves must have degree divisible by $2g-2$, after which we can use the same methods as \citep{mestrano} to deduce the original strong Franchetta conjecture.% We can for example immediately deduce that on a surface in $\bP^3$ of Picard number one and degree at least four, the general curve in any linear system on the surface is not hyperelliptic.

The hypothesis on the dimension of the locus of non-integral divisors might seem a little odd at first, but we will show that there are counterexamples to the theorem without this hypothesis. Moreover, in the case of very ample linear systems on a surface, the Castelnuovo-Kronecker theorem (see \citep{Castelnuovo}) implies that this hypothesis is satisfied except in the case where the image of the map from the surface to projective space is either ruled by lines or a projection of the Veronese surface.

Our proof strategy will be to first prove the analogue of the weak Franchetta conjecture, which is very easy in this setting, then to restrict our section to pencils, where we can use Tsen's theorem to show that any rational section must be a linear combination of the restriction of a line bundle on $X$ and the base locus of the pencil. We then show that in fact we can choose a line bundle on $X$ which gives rise to the section.

I would like to thank Dawei Chen and Joe Harris for their many helpful conversations on this topic. I would also like to thank Clifford Earle, Nicole Mestrano, Brendan Hassett, and Steve Kleiman for their help.

For us, a curve will be a connected projective scheme of dimension one. We work over $\mathbb{C}$, though the proof only needs an uncountable algebraically closed field of characteristic 0.

\section{Relative Picard Varieties}

For the basic facts about relative Picard varieties we recall in this section, we refer the reader to \citep{kleiman} except where otherwise noted. For an introduction to Brauer groups, see for example \citep{milne}.

Given a smooth projective morphism of varieties $\pi:\sD \to S$, we can define the relative Picard variety, which is a countable disjoint union of projective varieties $\Pic(\sD/S)^\eta$. The components of the relative Picard variety form a group, which we will call the relative Neron-Severi group, $\mathrm{NS}(\sD/S)$.

%\begin{cor}
%Let $\pi:\sD \to S$ be a smooth projective morphism with $S$ a smooth one-dimensional base. Let $L \in \Pic(\sD)$. If $L$ restricts to the trivial line bundle on the general fiber, then it restricts to the trivial line bundle on every fiber.
%\end{cor}
%\begin{proof}
%This follows from the valuative criterion for separatedness for the relative Picard functor.
%\end{proof}

\begin{prop} \label{relpic}
Let $\pi:\sD \to S$ be a smooth projective morphism. Let $L \in \Pic(\sD)$ restrict to the trivial line bundle on each fiber of $\pi$. Then $L \cong \pi^* L'$ for some $L' \in \Pic(S)$.
\end{prop}
\begin{proof}
Consider $L'=\pi_* L$. By Grauert's theorem (corollary (III, 12.9) of \citep{hartshorne})
, this is a line bundle, since there is one section of the trivial bundle on an integral scheme. Now consider \[\pi^* L'=\pi^* \pi_* L \to L\]. This is a nonzero map of invertible sheaves, and it is an isomorphism on fibers by Grauert's theorem, 
so it is an isomorphism of sheaves.
\end{proof}

The following is corollary 1.5 of \citep{mestranoramanan} in the case where the fibers are curves, but the proof is the same in general.

\begin{prop} \label{section_extend}
Let $\pi:\sD \to S$ be a smooth projective morphism. Let $\tau \in \mathrm{NS}(\sD/S)$. Let $\sigma:S \to \Pic^\tau(\sD/S)$ be a rational section of the natural map $\Pic^\tau(\sD/S) \to S$. Then $\sigma$ extends to a regular section.
\end{prop}

\begin{prop} \label{brauertorsion}
Let $\pi:\sD \to S$ be a smooth projective morphism with $S$ a smooth base. Let $\tau \in \mathrm{NS}(\sD/S)$. Let $\sigma:S \to \Pic^\tau(\sD/S)$ be a section. There is a natural number $m$ such that $\sigma^{\otimes m}$ comes from a line bundle on $\sD$.
\end{prop}
\begin{proof}
The obstruction to $\sigma$ coming from a line bundle on $\sD$ is an element of the Brauer group of $S$, and taking the tensor product of two sections adds these obstructions, but every element of the Brauer group of $S$ is torsion.
\end{proof}

\begin{prop} \label{tsen}
Let $\sD \to S$ a family of smooth projective varieties, with $S$ a smooth curve. Let $\tau \in \mathrm{NS}(\sD/S).$ Let $\sigma:S \to \Pic^\tau(\sD/S)$ be a section. Then there is a line bundle $L$ on $\sD$ which gives rise to $\sigma$.
\end{prop}
\begin{proof}
The obstruction to $\sigma$ coming from a line bundle on $\sD$ is an element of the Brauer group of $S$, but Tsen's theorem says that the Brauer group of a curve over an algebraically closed field is trivial.
\end{proof}

Suppose we have a family of projective varieties $\pi:\sD/S$ over an integral base such that the general fiber is smooth. By abuse of notation, we will refer to the relative Picard variety (resp.~Neron-Severi group) of the restriction of $\pi$ to the complement of the discriminant locus in $S$ as the relative Picard variety(resp.~Neron-Severi group) of $\pi$.

\section{Counterexamples}

In this section, we show that the hypothesis on the dimension of the locus of non-integral divisors in theorem \ref{mainresult} is necessary.

The first counterexample is very simple. Take the complete linear system of conics in $\bP^2$. There is certainly a rational section of the relative Picard variety which assigns to a smooth conic $C$ the line bundle $\cO_C(1)$, but there is no line bundle on $\bP^2$ which restricts to $\cO_C(1)$ on each conic.

Conics are somewhat exceptional, having genus 0, so we will rest easier once we have found a counterexample using curves of higher genus. Indeed, we will show that there are counterexamples with arbitrarily high genus.

Let $S$ be a very general double cover of $\bP^2$ branched over a sextic curve. Then $S$ is a K3 surface of Picard number 1, generated by the pullback of $\cO_{\bP^2}(1)$. Let $C$ be the preimage of a conic in $\bP^2$. Then $C$ is a hyperelliptic curve of genus 5, so $|C| \cong \bP^5$, since on a K3 surface, the dimension of a linear system of curves is equal to the genus. Since the linear system of conics in $\bP^2$ is five-dimensional, this means that every curve in the linear system is a double cover of a conic.

There is a rational section of $J^2 \to |C|$ which sends each curve to the line bundle on that curve giving rise to the double cover of $\bP^1$. On the other hand, it is easy to check that there is no line bundle on $S$ which has intersection number 4 with $C$, so this rational section cannot come from a line bundle on $S$.

By the Noether-Lefschetz theorem for weighted projective spaces, a double cover of $\bP^2$ branched along a very general curve of degree at least 6 has Picard number 1. Taking the preimage of a conic in such a surface will give a hyperelliptic curve. It is not difficult to show that all curves in the same linear system will again be double covers of conics. The same argument as above shows that the corresponding rational section of $J^2 \to |C|$ cannot come from a line bundle on the surface. Increasing the degree of the branch curve increases the genus of the curves in the linear system, so this provides us with counterexamples of arbitrarily high genus.

\section{Proof of the Main Theorem}

Let $X$ be a smooth projective variety, and $|D|$ a linear system. The universal divisor $\sD$ over $|D|$ maps to $X$, so we can pull back any line bundle $L$ on $X$ to the universal divisor. By the universal property of the relative Picard scheme, $L$ gives rise to a rational section of the relative Picard variety of line bundles, and its image is contained in some component $\Pic^\tau(\sD/|D|)$ with $\tau \in \mathrm{NS}(\sD/|D|)$.

We will let $|D|^{s}$ be the complement of the discriminant locus, and $\sD^s$ its preimage in $\sD$. By Bertini's theorem, if $|D|$ is basepoint-free, then $|D|^{s}$ is nonempty. For the rest of this section, we will assume that the hypotheses of theorem \ref{mainresult} hold.

We first note that the analogue of the weak Franchetta conjecture for basepoint-free linear systems is very easy.

\begin{lem}The relative Picard group of the universal divisor over a basepoint-free linear system on a projective variety $X$ is generated by $\Pic(X)$.
\end{lem}
\begin{proof}
Since $|D|$ is basepoint-free, the natural map $\sD \to X$ realizes $\sD$ as a projective bundle over $X$, so its Picard group is the direct sum of $\Pic(X)$ and the tautological quotient line bundle $\cO(1)$, but $\cO(1)$ is pulled back from $|D|$.\end{proof}

By proposition \ref{brauertorsion}, this means that there is an integer $m$ such that $\sigma^{\otimes m}$ comes from some line bundle on $X$, which we will call $L$.

%We now recall the following. % Try to find a reference

%\begin{lem}
%Let $\tilde{X}$ be the blowup of $X$ at a smooth subvariety. Then $\Pic(\tilde{X}) \cong \Pic(X) \oplus G$, where $G$ is the free Abelian group generated by the connected components of the exceptional locus.
%\end{lem}

\begin{prop}
For every $[D'] \in |D|^s$, we have $\sigma([D'])=[L_{D'}|_{D'}]$ for some $L_{D'} \in \Pic(S)$.
\end{prop}
\begin{proof}
%Consider a general pencil of divisors containing $D'$. By the hypothesis on the codimension of the non-integral locus, every divisor in the pencil is integral. By proposition \ref{section_extend}, we know that the rational section of the relative Picard variety is define on the entire pencil.

Fix $D' \in |D|$ and consider a general pencil containing $D'$. Let $\tilde{X}$ be the total space of the pencil, and $\bP^1$ the base. $\tilde{X}$ is the blowup of $X$ at the scheme-theoretic base locus of the pencil, which is smooth by a double application of Bertini's theorem, so \[\Pic(\tilde{X})\cong \Pic(X) \oplus \Z E_i\] where the $E_i$ are the connected components of the exceptional locus. Let $D_p$ denote the fiber over a point $p \in \bP^1$. The class of $D_p$ in $\Pic(\tilde{X})$ will be $D-E$, where $E=\sum E_i$ is the (reduced) exceptional divisor.

Let $C \subset \bP^1$ be the complement of the discriminant locus (or the complement of a point if the discriminant locus is empty). By proposition \ref{section_extend}, $\sigma$ is defined on all of $C$. Let $D_C$ be the preimage of $C$ in $\tilde{X}$. All the fibers of the map $\tilde{X} \to \bP^1$ are integral by hypothesis, so \[\Pic(D_C) \cong \Pic(\tilde{X})/(D-E) \cong \Pic(X) \oplus \Z E_i/(D-E)\] by the exact sequence of divisor class groups for an open subset ((II, 6.5) of \citep{hartshorne}).

By proposition \ref{tsen}, $\sigma$ comes from a line bundle on $D_C$, $\tilde{L}$. We can write \[\tilde{L} \equiv L'+\sum a_i E_i \pmod {D-E}\] with $L' \in \Pic(X)$ by the calculation of $\Pic(D_C)$. We know that \[mL'+m\sum  a_i E_i-L\] restricts to the trivial line bundle on the fibers of $D_C/C$, so by proposition \ref{relpic}, it is trivial in $\Pic(D_C)$ (since the class of a fiber is trivial), and hence a multiple of $D-E$ in $\Pic(\tilde{X})$ (say $n(D-E)$). We can rewrite this fact as the equation \[mL'-L+nD+\sum (m a_i-n) E_i=0\] in $\Pic(\tilde{X})$. Since all the $E_i$ are linearly independent from each other and from $\Pic(X)$ in $\Pic(\tilde{X})$, this means that in particular all the $a_i$ are equal, say to $a$. But $L'+ a \sum E_i=L'+aE$ has the same restriction to each fiber of $D_C/C$ as $L'+aD$, which is the pullback of a line bundle on $X$.
\end{proof}

Consider the set \[T=\{\tau' \in \mathrm{NS}(X):\tau'|_{|D|}=\tau\}\]. For each $\tau' \in T$, we have a restriction map \[\Pic^{\tau'}(X) \times |D|^s \to \Pic^\tau(\sD^s/|D|^s)\]. Each of these maps is proper over $|D|^s$ since \[\Pic^{\tau'}(X) \times |D|^s \to \Pic^\tau(\sD^s/|D|^s)\] is proper over $|D|^s$ (since $\Pic^{\tau'}(X)$ is proper) and the map \[\Pic^\tau(\sD^s/|D|^s) \to |D|^s\] is proper. In particular, each of the restriction maps has a closed image. By the above proposition, the image of $\sigma$ is contained in the union of these images. But $T$ is a countable set by Severi's theorem of the base, so by pulling back by $\sigma$  the images of the restriction maps for each $\tau' \in T$, we see that $|D|^s$ is a countable union of closed subvarieties, so one of them must be all of $|D|^s$, and hence there must be some $\tau'$ such that the image of $\sigma$ is contained in the image of $\Pic^{\tau'}(X)$.

Pick $\overline{L} \in \Pic^{\tau'}(X)$ with $\tau'$ chosen as above, and let $\sigma_{\overline{L}}$ be the corresponding section of $\Pic^\tau(\sD^s/|D|^s)$. By considering $\sigma-\sigma_{\overline{L}}$, we might as well assume that $\tau$ and $\tau'$ are both 0.

Now consider the map \[r:\Pic^0(X) \times |D|^s \to \Pic^0(\sD^s/|D|^s)\]. This is a morphism of abelian varieties over $|D|^s$ which preserves 0, so in particular, it's a group homomorphism. Let $K$ be the kernel of $r$. Let \[\pi:\Pic^0(X) \times |D|^s \to |D|^s\] be the projection onto the second factor. We will need the following lemma.

\begin{lem}
There is a nonempty open set $U \subset |D|^s$ such that $K \cap \pi^{-1}(x)$ is constant for $x \in U$.
\end{lem}
\begin{proof}
By generic flatness, there is an open set $V \subset |D|^s$ such that $K_x=K \cap \pi^{-1}(x) \subset \Pic^0(X)$ is a flat family of closed subvarieties. We will now restrict our attention to $V$.

Consider the component of the Hilbert scheme of closed subvarieties of $\Pic^0(X)$ which contains $K_x$. The tangent space to this point of the Hilbert scheme is given by $H^0(N_{K_x/\Pic^0(X)})$, but since $K_x$ is a closed subgroup, this normal bundle is a trivial bundle of rank equal to $c$, the codimension of $K_x$ in $\Pic^0(X)$. Therefore,  $h^0(N_{K_x/\Pic^0(X)})=cn$, where $n$ is the number of components of $K_x$.

We will now construct a flat family of embedded deformations of $K_x$ in $\Pic^0(X)$ such that its base dominates this component of the Hilbert scheme. Assume first that $K_x$ is connected. We note that we can identify the vector space $N_{K_x/\Pic^0(X),0}$ with $T_0(\Pic^0(X)/K_x)$. Let \[K'_x \subset \Pic^0(X) \times \Pic^0(X)\] be such that \[\pi_2^{-1}(\{a\})=K_x+a\], i.e.~$K_x$ translated by $a$. This is just the universal family of translates of $K_x$. There is an induced map from $\Pic^0(X)$, considered as the base of this family, to the Hilbert scheme of subschemes of $\Pic^0(X)$, and the differential of this map at 0 is given by the natural map \[T_0 \Pic^0(X) \to T_0(\Pic^0(X)/K_x) \cong H^0(N_{K_x/\Pic^0(X)})\], which is certainly surjective. Moreover, the kernel of this map consists of directions in which $a \in K_x$, or equivalently, directions in which $0 \in K_x+a$. In particular, any point near to $K_x$ but not equal to it, cannot be a subgroup, since it will not contain 0.

If $K_x$ is not connected, then nearby points of the Hilbert scheme will correspond to independent translations of each component of $K_x$, so again, none of the nontrivial deformations of $K_x$ can be a subgroup.
\end{proof}

Let $K_0$ be $K_x$ for the $x$ in the $U$ of the above lemma. We get a birational factorization \[\Pic^0(X) \times |D|^s \to \Pic^0(X)/K_0 \times |D|^s \dashedrightarrow \Pic^0(\sD^s/|D|^s)\] where the last arrow is a rational map which is birational onto its image. We know that $\sigma$ is contained in the closure of the image of this last map, and $\sigma$ is defined for all points of $|D|^s$ by proposition \ref{section_extend}, so we see that $\sigma$ factors birationally to give a map \[|D|^S \dashedrightarrow \Pic^0(X)/K_0.\] Since $|D|^s$ is an open subvariety of projective space and $\Pic^0(X)/K_0$ is an abelian variety, this map must be constant. We can therefore find an element $L'$ of $\Pic^0(X)$ such that $\sigma$ and $\sigma_{L'}$ agree on a dense open subset of $|D|^s$, and hence agree everywhere.

\bibliographystyle{plain}
\bibliography{divisor-franchetta-bibliography}

\end{document}